\documentclass[11pt,a4paper]{amsart}

\usepackage{amsmath,amscd}
\usepackage{amsthm}
\usepackage{amssymb}
\usepackage{graphics}
\usepackage{graphicx}
\usepackage{colortab} %Para tener colores
\usepackage{enumitem}
 
\usepackage{esint}
\usepackage{amscd}
\usepackage{pstricks}
\usepackage{pst-plot}
\usepackage{multido}
\usepackage{pst-coil}
\usepackage{colortab}

\usepackage[colorlinks=true,linkcolor=blue,citecolor=blue]{hyperref}
\usepackage{prettyref}
\theoremstyle{plain}
  \newtheorem{theorem}{\bf Theorem}[section]
  \newtheorem{proposition}[theorem]{\bf Proposition}
  \newtheorem{lemma}[theorem]{\bf Lemma}
  \newtheorem{corollary}[theorem]{\bf Corollary}

\theoremstyle{remark}
  \newtheorem{remark}[theorem]{\bf Remark}
  \newtheorem{definition}[theorem]{\bf Definition}
 
	 \newcommand{\R}{\mathbb{R}}

	  \newcommand{\esf}{\mathbb{S}}

\begin{document}

	%%%%%%%%%%%%%%%%%%%%%%%%%%%%%%%%%%%%%%%%%%%%%%%%%%%%%%%%%%%%%%%%%%%%%%%%%%%%%%%%%%%%%%%%%%
	\title{Characterization of Sobolev spaces on the sphere}
	%\author{J.A. Barcel\'{o}}\thanks{Departamento de Matem\'atica e Inform\'atica aplicadas a las Ingenier\'ias Civil y Naval, Universidad Polit\'ecnica de Madrid. E-mail: juanantonio.barcelo@upm.es} % \and  T. Luque\thanks{Departamento de An\'alisis Matem\'atico y Matem\'atica Aplicada, Universidad Complutense de Madrid. E-mail: t.luque@ucm.es.} \and M. C. Vilela\footnotemark[1]}
	
	\author{J. A. Barcel\'{o}}
	\address{Departamento de Matem\'atica e Inform\'atica aplicadas a las Ingenier\'ias Civil y Naval, Universidad Polit\'ecnica de Madrid.} 
	\email{juanantonio.barcelo@upm.es}

	\author{T. Luque}
	\address{Departamento de An\'alisis Matem\'atico y Matem\'atica Aplicada, Universidad Complutense de Madrid.} 
	\email{t.luque@ucm.es}
	\author{S. P\'erez-Esteva}
	\address{Instituto de Matem\'aticas, Unidad de Cuernavaca, Universidad Nacional Aut\'onoma de M\'exico.} 
	\email{spesteva@im.unam.mx }
	\keywords{Sobolev spaces, zonal Fourier multiplier, square function}
	\date{August 27, 2019.}
	\subjclass[2000]{Primary 46E35, 42B37; Secondary 42B35}
	\maketitle
	% 46E35 Sobolev spaces and other spaces of "smooth'' functions, embedding theorems, trace theorems
	% 	42B35  	Function spaces arising in harmonic analysis
	% 42B37 Harmonic analysis and PDE
	% 	Harmonic analysis in one variable 42A45  	Multipliers??
	%%%%%%%%%%%%%%%%%%%%%%%%%%%%%%%%%%%%%%%%%%%%%%%%%%%%%%%%%%%%%%%%%%%%%%%%%%%%%%%%%%%%%%%%%%%%
	
	\begin{abstract}We prove a characterization of the Sobolev spaces $H^\alpha$ on the unit sphere $\mathbb{S}^{d-1}$, where the smoothness index  $\alpha$ is any positive real number and $d\geq 2$. This characterization does not use differentiation and it is given in terms of  $([\alpha/2]+1)$-multidimensional square functions $S_\alpha$.  For  $[\alpha/2]=0,$ a function $f\in L^2(\mathbb{S}^{d-1})$ belongs to $H^\alpha(\mathbb{S}^{d-1})$ if and only if $S_\alpha (f)\in  L^2(\mathbb{S}^{d-1})$. If $n=[\alpha/2]>0$, the membership of $f$ is equivalent to the existence of $g_1,\cdots,g_n$ in $L^2(\mathbb{S}^{d-1})$ such that $S_\alpha(f,g_1,\ldots,g_n)\in L^2(\mathbb{S}^{d-1})$ and in this case,   $g_j=T_j((-\Delta_S)^j f)$, where  $T_j$ is a zonal Fourier multiplier in the sphere and $\Delta_S$ is the Laplace-Beltrami operator. The square functions $S_\alpha $  are based on averaging operators over euclidean balls (caps) in the sphere that may be viewed as zonal multipliers.  The results in the paper are in the spirit of the characterization of fractional Sobolev spaces given in $\mathbb{R}^d$ proved in \cite{AMV}. The development of the theory is fully based on zonal Fourier multipliers and  special functions. 
	\end{abstract}

\section{ Introduction}
The interest in characterizing the Sobolev spaces $W^{\alpha,p}(\mathbb{R}^d)$ without involving distributional derivatives goes back to the 1960s, with the early works of Stein \cite{S}, Strichartz \cite{Strich} and the generalization of Bagby \cite{Bagby}. These first characterizations are in terms of certain mix $L^{p,2}$ norms of a first difference, with $1<p<\infty$ and $\alpha>0$. In \cite[p.139]{Sbook}, Stein again presents an alternative characterization of Sobolev spaces $W^{1,p}(\mathbb{R}^d)$, $1<p<\infty$, in terms of the order of smallness of the  $L^p$-modulus of continuity, which can be extended to $p=2$ and $0<\alpha<2$. Later on, Dorronsoro \cite{Dorr} gave also a characterization by using a mixed norm estimate of a mean oscillation that is defined as the difference of the function from a certain polynomial approximation. This characterization is for $\alpha>0$ and $1<p<\infty$.

As it is pointed out in \cite[p.243]{David}, the characterizations mentioned above can be read in terms of square functions with different forms, so that differentiability properties of functions in $\mathbb{R}^n$ are characterized in terms of the boundedness of these square functions in certain spaces. This idea is also seen in \cite{Stein_square}.

More recently, the interest of defining and studying  Sobolev type spaces on  metric measure spaces, where there is a lack of a differential structure, has led to new non-differential characterizations of the classical Sobolev spaces.  See for example the works of Hajłasz \cite{Hajl1,Hajl2}, which consider the concepts of upper gradients and weak upper gradients. \cite{Shan2} and \cite{HST} also follow this approach. 

Again in the spirit of the characterization via square functions and also in the context of metric spaces, we refer the work of Alabern, Mateu and Verdera \cite{AMV}. The authors gave a characterization of the Sobolev spaces $W^{\alpha,p}(\mathbb{R}^d) $ of any positive order of smoothess $\alpha$, using a square function associated with this $\alpha$. More precisely, given a locally integrable function $f$ they define  a square function  $S_\alpha (f)$, for $\alpha\in(0,2)$, and a quadratic multiscale operator $S_\alpha(f,g_1,\cdots,g_n)$, for $\alpha\geq 2$ and where  $n=[\alpha/2]$. Then they prove that a function $f$ in $ L^{p}(\mathbb{R}^d)$, $1<p<\infty$, belongs to the Sobolev space $ W^{\alpha,p}(\mathbb{R}^d)$ if and only if there exist $g_1,\cdots,g_n\in L^{p}(\mathbb{R}^d)$ such that $S_\alpha(f,g_1,\cdots,g_n)\in L^{p}(\mathbb{R}^d)$ (resp $S_\alpha (f)\in L^{p}(\mathbb{R}^d)$ if $\alpha\in (0,2)$). 
The square functions defined in \cite{AMV} are given by
\begin{equation}\label{AMV1}
S_\alpha f(x)^2= \int_0^\infty \left|\fint_{B(x,t)}\frac{f(y)-f(x)}{t^\alpha}\,dy\right|^2\frac{dt}{t},\quad \text{if }\alpha\in (0,2),
\end{equation}
where $B(x,t)$ is the ball of center $x\in\mathbb{R}^d$  and radius $t>0$, and where the barred integral on the set $B(x,t)$ stands for the mean over that set. If $ \alpha\in(2n,2(n+1)), \; n=1, 2, \dots$, the square function is given by
\begin{align}\label{splana}
&S_\alpha (f,g_1,\cdots,g_n)(x)^2 \nonumber\\
=&\int_0^\infty \left|  \fint_{B(x,t)} \frac{f(y)-f(x)-\sum_{i=1}^{n}g_i(x)\vert y-x\vert^{2i}}{t^\alpha}\,dy\right|^2\frac{dt}{t},
\end{align}
In the particular case $\alpha = 2n, \; n=1,2, \dots$,  $g_n$ in the right hand side of \eqref{splana} is replaced by the more regular function $\fint_{B(x,t)}g_n(y)dy.$ We remark that these square functions follow the same lines that the ones presented at the begining of this introduction; in particular, they define a mean oscillation that is given as the difference of the function from its Taylor approximation.

Since these definitions depend only on the metric of $\mathbb{R}^d$ and the Lebesgue measure, this characterization can be adapted to define Sobolev spaces of any metric measure space. New developments and extensions of this characterization using different operators of the Littlewood-Paley theory can be found in \cite{DLYY}, \cite{HYY1} and \cite{HYY2}.

The purpose of this paper is to give characterizations of the Sobolev spaces in the unit sphere, $\mathbb{S}^{d-1}$, via  square functions. We will prove that the natural versions in the  sphere of the square functions defined in \eqref{AMV1} and \eqref{splana} will provide such characterization. In fact , in \cite[Final Remarks]{AMV} the authors point out as an open problem the extension of their results to manifolds and the present work is a step in this direction. 

Our real motivation was the characterization of solutions in $\mathbb{R}^d$ of the Helmholtz equation $\Delta u+u=0$ that arise  as  the Fourier extension operator  of distributions in the sphere. The classic example of this is the space of Herglotz wave functions, that are  the image under the Fourier extension operator of functions in $L^2(\mathbb{S}^{d-1})$. These functions play an important role in scattering theory (see \cite{CK}) and can be characterized by the  classical condition of Hartman and Wilcox \cite{HW}  involving  the  growth condition
\begin{equation}\label{HW}
\limsup_{R\rightarrow\infty}\,\frac{1}{R}
\int_{|x|<R}|u(x)|^2\,dx
<\infty.
\end{equation} 
In \cite{preprint},we study the Herglotz wave functions that are  the image under the Fourier extension operator of functions in the Sobolev space $W^{\alpha, 2}(\mathbb{S}^{d-1})$, $\alpha>0$, and using the characterization of these spaces obtained in here, we are able to deduce a condition in the spirit of \eqref{HW}. 

Although the characterization  of $W^{\alpha, 2}(\mathbb{S}^{d-1})$ presented here  extends naturally  the one given in \cite{AMV} for $\mathbb{R}^d$, the methods are completely different. The reader will see that the results in the spherical setting was accomplished by a careful study of certain special functions and the subtle use of classical analysis tools.

In order to present this characterization, we start with the definition of the Sobolev spaces $W^{\alpha, 2}(\mathbb{S}^{d-1})$, that for simplicity will be denoted by $H^{\alpha}(\mathbb{S}^{d-1})$. Throughout this paper $L^2(\esf^{d-1})$ will stand for the space of square integrable  functions in the sphere  provided with the uniform measure  $d\sigma$. Let   $\Delta_S$ denote the Laplace-Beltrami operator on the unit sphere. Then the Sobolev space  $H^\alpha (\esf^{d-1})$, $\alpha>0$, is defined as  the space of all functions $f\in L^2(\esf^{d-1})$ such that $(-\Delta_S)^{\alpha/2} f\in  L^2(\esf^{d-1})$.  We provide   $H^\alpha (\esf^{d-1})$ with the Hilbert  space norm
\begin{equation*}
 \|f\|_{H^{\alpha}(\mathbb{S}^{d-1})}^2=\|f\|_{L^2(\esf^{d-1})}^2+\|(-\Delta_S)^{\alpha/2} f\|_{L^2(\esf^{d-1})}^2.
 \end{equation*} 
 
The space of all spherical harmonics of degree $\ell$  is denoted by $\mathbb{H}_\ell^d$. In addition, $\{Y_\ell^j :\ell=0,1,\ldots,\, 1\leq j\leq \nu(\ell)\}$ will denote a orthonormal basis of
spherical harmonics for $L^2(\esf^{d-1})$, where $\nu(\ell)$ is the dimension of the space $\mathbb{H}_\ell^d$. It is well known that each $Y_\ell^j$ is an eigenfunction of   $(-\Delta_S$) corresponding to the eigenvalue $\ell(\ell+d-2)$. Then we can rewrite $H^\alpha(\mathbb{S}^{d-1})$ as those functions $f\in L^2(\mathbb{S}^{d-1})$ such that
\begin{equation}\label{def:H}
\|f\|_{H^{\alpha}(\mathbb{S}^{d-1})}^2=\sum_{\ell=0}^\infty\sum_{j=1}^{\nu(\ell)} (1+\ell^{\frac{1}{2}}(\ell+d-2)^{\frac{1}{2}})^{2\alpha}|\hat f_{\ell j}|^2<\infty,
\end{equation}
where $f$ always admits the following representation
  \begin{equation}\label{representationFunct}
f=\sum_{\ell=0}^\infty\sum_{j=1}^{\nu(\ell)} \hat{f}_{\ell j}\,Y_\ell^j, \quad \hat{f}_{\ell j}:=\int_{\mathbb{S}^{d-1}}f(\xi)\,\overline{Y_\ell^j}(\xi)\,d\sigma(\xi).
\end{equation}

We first consider the space $H^\alpha(\mathbb{S}^{d-1})$,  for $0<\alpha<2$, and we define the square function similar to \eqref{AMV1}; that is, for $f$ integrable function on $\mathbb{S}^{d-1}$
\begin{equation}\label{def:Sfunction}
S_\alpha(f)^2(\xi):=\int_0^\pi\left|\frac{A_tf(\xi)-f(\xi)}{t^\alpha} \right|^2\frac{dt}{t}, 
\end{equation}
where
\begin{equation}\label{def:At}
A_tf(\xi):=\fint_{C(\xi,t)}f(\tau)\,d\sigma(\tau)
\end{equation}
denotes the mean of $f$ on the spherical cap centred at $\xi\in\mathbb{S}^{d-1}$ and angle $t\in(0,\pi)$ given by
\begin{equation}\label{def:casquete}
C(\xi,t):=\{\eta\in\mathbb{S}^{d-1}:\xi\cdot\eta\geq\cos t\}.
\end{equation}

The characterization obtained is the following.
	\begin{theorem}\label{prop:rango02}
	Let $0<\alpha<2$ and $f\in L^2(\mathbb{S}^{d-1})$. Then $f\in H^\alpha(\mathbb{S}^{d-1})$ if and only if $S_\alpha f\in L^2(\mathbb{S}^{d-1})$. 
	\end{theorem}

%To extend Proposition \ref{prop:rango02} to $H^2(\mathbb{S}^{d-1})$ and higher orders of %smoothness $\alpha$, we have find a replacement when  $\alpha\geq 2$ for $f(\xi)-f(\eta) $ in the %definition of $S_\alpha$. To get a hint do this we argue as  in the lines of  \cite{AMV}: let $f:\mathbb{S}%\rightarrow\C$ is smooth and $\xi\in \mathbb{S}^{d-1}$, for simplicity assume that $\xi=e_d$. Take %the chart of $\mathbb{S}^{d-1}$ at $e_d$ given by $\varphi (y)=(y,\sqrt{1-\vert y\vert^2}).$
%Let $P$ the second degree  homogeneous part of the Taylor expansion of $f\circ\varphi$ at $y=0$. %Then we can write
%\begin{equation*}
%P(y)=H(y)+ \frac{\Delta (f\circ\varphi) (0)}{2(d-1)}\vert y\vert^2
%\end{equation*}
%where $H$ is harmonic.
%Thus
%\begin{equation*}
%A_t (P\circ\varphi^{-1})(e_d)=\frac{\Delta (f\circ\varphi) (0)}{2(d-1)}A_t (\vert \varphi^{-1})\vert^2)%%(e_d).
%\end{equation*}
%The choice of $\varphi$ is such that  $\Delta_S f(e_d)= \Delta (f\circ\varphi) (0)$ and $A_t (\vert %\varphi^{-1}\vert^2)(e_d)=A_t(\vert e_d-\cdot\vert^2)+O(t^4).$
To characterize the elements of $H^\alpha \left(  \mathbb{S}^{d-1} \right)$ for $\alpha\geq 2$, we first need to extend the definition of the square function.
\begin{definition}\label{def:square} Let $\alpha>0$. The square function is defined 
by 
 \eqref{def:Sfunction} whenever $\alpha<2$. The case $\alpha\geq 2$ is extended as follows: let $n\in  \mathbb{N}$ such that  $2n\leq\alpha<2(n+1)$ and $f,g_1,\ldots,g_n$ integrable functions 
on $\mathbb{S}^{d-1}$. Then the square function is given by
\begin{align*}%\label{def:Sgeneral}
S_\alpha(f,g_1,\ldots,g_n)(\xi)^2&:=\int_0^\pi\bigg|A_tf(\xi)-f(\xi)-\sum_{k=1}^ng_k(\xi)A_t \big( \left|  \xi - \cdot \right|^{2k}    \big)(\xi)   \bigg|^2 \frac{dt }{t^{2 \alpha +1}},
\end{align*}
for $2n<\alpha<2(n+1)$ and  
\begin{align*}
S_{2n}(f,g_1,\ldots,g_n)(\xi)^2&:=\int_0^\pi\bigg|A_tf(\xi)-f(\xi)-\sum_{k=1}^{n-1}\,g_k(\xi) A_t \big( \left|  \xi - \cdot \right|^{2k}    \big)(\xi)\nonumber\\
&-A_t(g_n)(\xi)A_t\big( \left|  \xi - \cdot \right|^{2n}    \big)(\xi)   \bigg|^2 \frac{dt }{t^{2 \alpha +1}},%\label{def:Spar}
\end{align*}
$\xi\in\mathbb{S}^{d-1}.$ 
\end{definition}

 We have the following characterization of $H^\alpha \left(  \mathbb{S}^{d-1}\right)$, for $\alpha \geq 2$.

\begin{theorem}\label{pepe} 	Let $\alpha \geq 2$ and $n\in \mathbb{N}$ such that  $2n \leq \alpha < 2(n+1)$.  Then   $f \in H^\alpha \left(  \mathbb{S}^{d-1} \right)$ if and only if there exist  $g_1,g_2, \cdot \cdot \cdot, g_n \in L^2 \left( \mathbb{S}^{d-1}  \right) $ such that $f,S_\alpha (f,g_1, g_2, \cdot \cdot \cdot, g_n ) \in L^2 \left( \mathbb{S}^{d-1}  \right) $, where $S_\alpha$ is the square function introduced in Definition \ref{def:square}.
	\end{theorem}

	The proof of Theorem \ref{pepe} will be based on the following proposition.
\begin{proposition}\label{teorema1}
	Let $\alpha>0$ and $n$ the non-negative integer number such that  $2n\leq\alpha<2(n+1)$. Let $S_\alpha$ be  the square function introduced in Definition \ref{def:square}. If $n=0$, then
		\[\|f\|_{H^{\alpha}(\mathbb{S}^{d-1})}\sim_{d,n}\|S_\alpha (f)\|_{L^2(\mathbb{S}^{d-1})}.\]
	 If $n\geq 1$, then for $k=1,2\ldots n$, there exist isomorphisms $T_k$ of $L^2(\mathbb{S}^{d-1})$ such that for $f\in H^{2n}(\mathbb{S}^{d-1})$
	\[\|f\|_{H^{\alpha}(\mathbb{S}^{d-1})}\sim_{d,n}\|S_\alpha (f,T_1((-\Delta_S)f),\ldots ,T_n((-\Delta_S)^nf))\|_{L^2(\mathbb{S}^{d-1})}.\]
 Moreover, each  $T_k$ is a zonal Fourier multiplier of $L^2(\mathbb{S}^{d-1})$ with multiplier  given in \eqref{betamult} (see definition \eqref{zonal_fourier} below).
\end{proposition}
Observe that Theorem \ref{prop:rango02} particularises Proposition \ref{teorema1} for $n=0$. 
\paragraph{\emph{Outline}.} The general organization of this paper is as follows. The second section is devoted to the proof of Proposition \ref{teorema1} and Theorem \ref{pepe}. The proof of these resuts need certain technical auxiliary lemmas that are proved in the last section.
\paragraph{\emph{Notation}.} Given a function $f$ we use Lagrange's notation to denote the higher order derivatives; thus $f^{(n)}$ defines the nth derivative. For non-negative quantities $a$ and $b$, we write $a  \lesssim b$ ($a\gtrsim b$) if $a\leq c b$ ($a\geq c b$) for some positive numerical constant $c>0$. We write $a\sim b$ if both $a   \lesssim  b$ and $a \gtrsim  b$ hold. In order to indicate the dependence of the constant $c$ on some parameter $n$ (say), we write $a\lesssim_n b$, $a\gtrsim_n b$ or $a\sim_n b$.
\section{Proofs}
\subsection{Proof of Proposition \ref{teorema1}}

  The proof of  Proposition \ref{teorema1} and Theorem \ref{prop:rango02} is based on two essential points. The first is the fact that the average operator $A_t$ involved on the square function is a zonal Fourier multiplier of $L^{2}(\mathbb{S}^{d-1})$; more precisely, given a bounded  sequence $\{m_\ell\}_{\ell =0,1, \cdots}$, a zonal Fourier multiplier associated to $m$ in this setting is the operator T defined as
 \begin{align}\label{zonal_fourier}
Tf(\xi)=\sum_{\ell=0}^\infty\sum_{j=1}^{\nu(\ell)}m_\ell\, \hat{f}_{\ell j}\,Y_\ell^j(\xi),
\end{align}
for  $f\in L^2(\esf^{d-1})$. The second reduces the proof of these results to a problem of multipliers of $L^{2}(\mathbb{S}^{d-1})$, that is indeed solved in the technical Lemmas \ref{prop:estima} and \ref{prop:estimaalpha2n} given in the next section.

To present the first point of the proof, we introduce the Legendre polynomials of degree $\ell$ in $d$ dimensions that are denote by $P_{\ell,d}$. These polynomials are proportional to the Jacobi polynomials $P_\ell^{\alpha,\beta}$; in fact, we have that
\begin{equation}\label{def:Legendre}
P_{\ell,d}=\frac{\ell!\, \Gamma(\frac{d-1}{2})}{\Gamma(\ell+\frac{d-1}{2})}P_{\ell}^{\frac{d-3}{2},\frac{d-3}{2}}
\end{equation}
where $\Gamma$ denotes the gamma function. See \cite[pp.39]{AH}) for more details on their definition.

\begin{lemma}\label{lemma:At}The operator $A_t$, $t \in (0,\pi)$, is a zonal Fourier multiplier with multiplier
	\begin{equation}\label{mutliplierAt}
	m_{\ell,t}:=C_{t,d}\int_{\cos t}^1P_{\ell,d}(s)(1-s^2)^{\frac{d-3}{2}}\,ds,\;\;\; \ell= 0,1, \cdots,
	\end{equation}
	where $C_{t,d}:=\frac{|\mathbb{S}^{d-2}|}{|C(\xi,t)|}$ and  $C(\xi,t)$ is the spherical cap defined in \eqref{def:casquete}.
	\end{lemma}
	\begin{corollary}\label{corolarioAt} 
	The following equalities and estimates hold.
	\begin{equation}\label{medida_casquete}
	|C(\xi,t)|=|\mathbb{S}^{d-2}|\,\int_{\cos t}^1(1-s^2)^{\frac{d-3}{2}}\,ds,
	\end{equation}
	where $|C(\xi,t)|$ is the measure of $C(\xi,t)$;
	\begin{equation}\label{Cte}
C_{t,d} \sim_d t^{1-d};
	\end{equation}
	\begin{equation}\label{Atlap}
	A_t((-\Delta_S)^kY_\ell)(\xi)=\ell^k(\ell+d-2)^k\, m_{\ell,t} Y_\ell(\xi),\quad k=1,2,\ldots, \;\;  
	\end{equation}
where $Y_\ell\in \mathbb{H}_\ell^d$ and $m_{\ell,t}$ is the multiplier given by \eqref{mutliplierAt};
	\begin{equation}\label{Atfunct}
A_t(|\xi-\cdot|^{2k})(\xi)=2^kC_{t,d}\int_{\cos t}^1(1-s)^k(1-s^2)^{\frac{d-3}{2}}\,ds,
	\end{equation}
	$$ \backsim \frac{2^k (d-1)}{2k+d-1}(1-\cos t)^k,\quad k=1,2,\ldots.$$
	
\end{corollary}
\emph{Proof of Proposition \ref{teorema1}.} 
Let $\alpha>0$ and $n$ the non-negative integer number such that  $2n\leq\alpha<2(n+1)$. We focus in the case $n\geq 1$, since the case $n=0$ is a particular situation that will be clarified later. We first define the zonal Fourier multipliers
\begin{equation*}
T_k(f)=\sum_{\ell=0}^\infty\sum_{j=1}^{\nu(\ell)}\beta_{k,\ell} \hat{f}_{\ell,j}Y_\ell^j,\quad k=1,2,\ldots,n.
	\end{equation*}
	where
	\begin{equation}\label{betamult}
			\beta_{k,\ell}=\frac{(-1)^kP_{\ell,d}^{(k)}(1)}{k!2^k\ell^k(\ell+d-2)^k}.
			\end{equation}
Since (see \cite[pp.58]{AH})
\begin{equation}\label{cotaderaiva_0}
P_{\ell,d}^{(k)}(1)=\frac{\ell!(\ell+k+d-3)!\Gamma(\frac{d-1}{2})  }{2^k(\ell-k)!(\ell+d-3)!\Gamma(k+\frac{d-1}{2}) }\sim\ell^k(\ell+d-2)^k,
\end{equation}

we have that for $k$ fixed $\beta_{k,\ell}\sim_{d,n} 1$ and hence, $T_k$ and $T_k^{-1}$ are continuous in   $L^2(\mathbb{S}^{d-1}).$

We distinguish two cases:  $\alpha\in(2n, 2(n+1))$ and $\alpha=2n$. In the first case, since $T_k$ is a zonal multiplier and $f$ admits representation \eqref{representationFunct}, by  \eqref{mutliplierAt}, \eqref{medida_casquete}, \eqref{Atfunct},
\begin{equation}\label{agosto}
P_{\ell,d}(1)=1=C_{t,d}\int_{\cos t}^1 (1-s^2)^{\frac{d-3}{2}}ds, \;\;\; \ell= 0,1, \ldots ,
\end{equation}
and 
\begin{equation}\label{laplaciano}
(-\Delta_S)^jY_\ell(\xi):=\ell^j(\ell+d-2)^j Y_\ell(\xi), \; \; \; j=1,2,\ldots,\quad Y_\ell\in\mathbb{H}_\ell^d,
\end{equation}
we have that  
\begin{align}
A_tf(\xi)&-f(\xi)-\sum_{k=1}^n\, T_k ((-\Delta_S)^k f)(\xi) A_t \big( \left|  \xi - \cdot \right|^{2k}    \big)(\xi)=\sum_{\ell=1}^\infty\sum_{j=1}^{\nu(\ell)}(m_{\ell,t}-1)\,\hat{f}_{\ell j}\,Y_\ell^j(\xi)\nonumber\\
-&\sum_{\ell=1}^\infty\sum_{j=1}^{\nu(\ell)}\sum_{k=1}^n\bigg[\beta_{k,\ell} \ell^k(\ell+d-2)^k\frac{2^k|\mathbb{S}^{d-2}|}{|C(\xi,t)|}\int_{\cos t}^1(1-s)^k(1-s^2)^{\frac{d-3}{2}}\,ds\bigg]\,\hat{f}_{\ell j}\,Y_\ell^j(\xi)\nonumber\\
=&\sum_{\ell=1}^\infty\sum_{j=1}^{\nu(\ell)} M_{\ell,t}\,\hat{f}_{\ell j}\,Y_\ell^j(\xi)\label{cason0},
\end{align}
where 
\begin{equation}\label{multiplicador1}
	M_{\ell,t}=C_{t,d}\int_{\cos t}^1\bigg[P_{\ell,d}(s)-P_{\ell,d}(1)-\sum_{k=1}^{n}c_{k,\ell}(1-s)^{k}\bigg](1-s^2)^{\frac{d-3}{2}}\,ds,
	\end{equation}
and
\begin{equation}\label{eq:coeficientes}
c_{k,\ell}:=\frac{(-1)^kP_{\ell,d}^{(k)}(1)}{k!},\quad k=1,2,\ldots,n,
	\end{equation}
	Notice that $P_{\ell,d}(s)-P_{\ell,d}(1)-\sum_{k=1}^{n}c_{k,\ell}(1-s)^{k}$ is the residue of the Taylor approximation  of $P_{\ell,d}$ of order $n$ at $s=1$. 
	
	Observe also that the case $n=0$ admits exactly the same expression that in \eqref{cason0} but removing the sum; that is,
\[A_tf(\xi)-f(\xi)=\sum_{\ell=1}^\infty\sum_{j=1}^{\nu(\ell)}M_{\ell,t}\,\hat{f}_{\ell j}\,Y_\ell^j(\xi),\]
with 
\[M_{\ell,t}=C_{t,d}\int_{\cos t}^1\bigg[P_{\ell,d}(s)-P_{\ell,d}(1)\bigg](1-s^2)^{\frac{d-3}{2}}\,ds.\]
It follows then
\begin{align}
& \|S_\alpha f\|_{L^2(\mathbb{S}^{d-1})}^2 =\int_{\mathbb{S}^{d-1}}|S_\alpha (f)(\xi)|^2\,d\sigma(\xi) \nonumber\\ & =\int_{\mathbb{S}^{d-1}}\int_0^\pi\bigg|\sum_{\ell=1}^\infty\sum_{j=1}^{\nu(\ell)} M_{\ell,t}\,\hat{f}_{\ell j}\,Y_\ell^j(\xi)\bigg|^2\frac{dt}{t^{2\alpha+1}}d\sigma(\xi)\nonumber\\
&=\sum_{\ell=1}^\infty\sum_{j=1}^{\nu(\ell)}|\hat{f}_{\ell j}|^2\int_0^\pi |M_{\ell,t}|^2\frac{dt}{t^{2\alpha+1}},\label{inter_Fubini}
\end{align}
Therefore, Proposition \ref{teorema1} in this case follows whenever 
\begin{equation*}
\int_0^\pi |M_{\ell,t}|^2\frac{dt}{t^{2\alpha+1}}\sim_{d,n}\ell^{2\alpha}
\end{equation*}  and this is proved in Lemma \ref{prop:estima}  below.

The case $\alpha=2n$, $n=1,2, \ldots$ repeats the same argument, but picking the second square function introduced in Definition \ref{def:square} and  property \eqref{Atlap} in Lemma \ref{lemma:At} is now needed. 
It can easily be seen that
\begin{align*}
&A_tf(\xi)-f(\xi)-\sum_{k=1}^{n-1}\, T_k (-\Delta_S)^k f)(\xi) A_t \big( \left|  \xi - \cdot \right|^{2k}    \big)(\xi) \\ &-T_n  ( A_t \left( \left(-\Delta_S  \right)^n f \right))(\xi) A_t \big( \left|  \xi - \cdot \right|^{2n}    \big)(\xi) 
=\sum_{\ell=1}^\infty\sum_{j=1}^{\nu(\ell)} N_{\ell,t}\,\hat{f}_{\ell j}\,Y_\ell^j(\xi),%\label{intteor1},
\end{align*}
where
\begin{align}\label{multiplicador2}
N_{\ell,t}&=C_{t,d}\int_{\cos t}^1\bigg[P_{\ell,d}(s)-P_{\ell,d}(1)-\sum_{k=1}^{n-1}c_{k,\ell}(1-s)^{k}\bigg](1-s^2)^{\frac{d-3}{2}}\,ds\nonumber\\
&-c_{n,\ell}\,C_{t,d}^2\int_{\cos t}^1P_{\ell,d}(s)(1-s^2)^{\frac{d-3}{2}}\,ds\int_{\cos t}^1(1-r)^n(1-r^2)^{\frac{d-3}{2}}\,dr,
\end{align}
and $c_{k,\ell}$ are given by (\ref{eq:coeficientes}). Observe that if we compare to the previous case, $N_{\ell,t}$ involves a finer approximation of $P_{\ell,d}$ than its  Taylor expansion of order $n$ at $s=1$.  Based on  Lemma \ref{prop:estimaalpha2n} proved below we have that
	\begin{equation*}
\int_0^\pi |N_{\ell,t}|^2\frac{dt}{t^{2\alpha+1}}\sim_{d,n}\ell^{2n}
\end{equation*} 
	and the proof of Proposition \ref{teorema1} is now complete.
\hfill$\Box$

We state now the key results that complete the proof of Proposition \ref{teorema1}.

\begin{lemma}\label{prop:estima}
	Let $t\in (0,\pi)$. The sequence $\left\{M_{\ell,t}\right\}_{\ell=1,2, \ldots}$  defined by (\ref{multiplicador1})
	%\begin{equation}\label{multiplicador1}
	%M_{\ell,t}=C_{t,d}\int_{\cos t}^1\bigg[P_{\ell,d}(s)-P_{\ell,d}(1)-\sum_{k=1}^{n}b_k(1-s)^{k}\bigg](1-s^2)^{\frac{d-3}{2}}\,ds, \;\; \ell \geq 1,
	%\end{equation}
	%and $M_{0,t}=0$, 
	%with $C_{t,d}$ defined in \eqref{Cte} and
	%\begin{equation}\label{eq:coeficientes}
%b_k:=\frac{(-1)^kP_{\ell,d}^{(k)}(1)}{k!},\quad k=1,2,\ldots,n,
	%\end{equation}
	verifies the following estimate
	\begin{equation}\label{clave2}
	I_{\alpha,n}(\ell)=\int_0^{\pi}  |M_{\ell,t}|^2\,\frac{dt}{t^{2\alpha+1}}\sim_{d,n} \ell^{2\alpha}
	\end{equation}
	whenever $2n<\alpha<2(n+1)$ and $n=0,1, \ldots$.
\end{lemma}

\begin{lemma}\label{prop:estimaalpha2n}
	Let $t\in(0,\pi)$. The sequence $\left\{N_{\ell,t}\right\}_{\ell=1,2, \ldots}$  defined by (\ref{multiplicador2})
	%\begin{align}\label{multiplicador2}
	%N_{\ell,t}&=C_{t,d}\int_{\cos t}^1\bigg[P_{\ell,d}(s)-P_{\ell,d}(1)-\sum_{k=1}^{n-1}b_k(1-s)^{k}\bigg](1-s^2)^{\frac{d-3}{2}}\,ds\nonumber\\
	%&-b_nC_{t,d}^2\int_{\cos t}^1P_{\ell,d}(s)(1-s^2)^{\frac{d-3}{2}}\,ds\int_{\cos t}^1(1-r)^n(1-r^2)^{\frac{d-3}{2}}\,dr,
	%\end{align}
	%with $C_{t,d}$ introduced in \eqref{Cte} and coefficients $b_k$ defined in \eqref{eq:coeficientes}, 
	verifies the following estimate
	\begin{equation}\label{clave3}
	J_{n}(\ell)=\int_0^{\pi}  |N_{\ell,t}|^2\,\frac{dt}{t^{4n+1}}\sim_{d,n} \ell^{4n}, \;\;\; n=1,2, \ldots.
	\end{equation}
\end{lemma}
\begin{remark}
	Observe that $m_{\ell,t}$, $M_{\ell,t}$ and $N_{\ell,t}$ do not depend on $\xi\in\mathbb{S}^{d-1}$ since $C_{t,d}$ defined in \eqref{Cte} is indeed independent of $\xi$. This is detailed in \eqref{medida_casquete}.
\end{remark}
We postpone the technical  proofs of Lemma \ref{lemma:At}, Corollary \ref{corolarioAt},     Lemmas \ref{prop:estima} and \ref{prop:estimaalpha2n} to  the next section.

\subsection{ Proof of Theorem \ref{pepe}.} The proof of Theorem \ref{pepe} precises a regularization of the operators $S_\alpha$ so we can apply Proposition \ref{teorema1}. This can be achieved smoothing $f, g_1\cdots g_n$ in $L^2(\mathbb{S}^{d-1})$. To do this, for convenience,  we choose the  Poisson transform as an approximation of the identity
\begin{equation}\label{transformada_Poisso}
P_rf(\xi)=\int_{\mathbb{S}^{d-1}}p_r(\eta, \xi) f(\eta) d \sigma (\eta),
\end{equation}
where $p_r(\eta , \xi)$ is the Poisson kernel in the unit ball or $\R^{d}$
$$p_r(\eta , \xi)= \frac{1-r^2}{|\mathbb{S}^{d-1}|\left|  r \xi - \eta \right|^{d}},$$ 
$0<r<1$ and $\xi,  \eta \in \mathbb{S}^{d-1}$.  

Let $ f \in L^2\left( \mathbb{S}^{d-1}  \right)$ and $S_\alpha$ the square function given in Definition \ref{def:square}, then the following properties hold (see \cite{Axler}).
\begin{itemize}
	\item[i)] $\Vert P_r f\Vert_{L^2(\mathbb{S}^{d-1})} \lesssim    \Vert f\Vert_{L^2(\mathbb{S}^{d-1})}.$ \vspace{0.2cm}
	\item[ii)] $P_r f\rightarrow f \text{ when $r \rightarrow 1^{-}$ in }  L^2(\mathbb{S}^{d-1}) \, \text{ and almost everywhere}.$ \vspace{0.2cm}
	\item[iii)] 
	$P_r f\in C^{\infty}(\mathbb{S}^{d-1})\subset H^\alpha(\mathbb{S}^{d-1}),$ with continous inclusion. In particular, $S_\alpha (P_rf,(-\Delta_S) P_r f\cdots,(-\Delta_S)^n P_r f)(\xi)$, $\xi\in\mathbb{S}^{d-1}$ is well defined. \vspace{0.1cm}
	\item[iv)] $p_r(\eta,\xi)=\sum_{\ell =0}^\infty \sum_{j=1}^{\nu (\ell)} r^\ell \overline{  Y_\ell^j(\eta)}Y_\ell^j(\xi)$,  which implies that $P_r$ is a zonal Fourier multiplier with multiplier  $\lbrace r^\ell\rbrace_{\ell=0,1, \ldots}$. \vspace{0.2cm}
	\item[v)] If $ g_1, \cdot \cdot \cdot, g_n, \in L^2 \left( \mathbb{S}^{d-1} \right)$ then
\begin{align}\label{propclave}
S_\alpha \left(P_rf, P_rg_1,\ldots P_rg_n   \right) (\xi) \lesssim  P_rS_\alpha \left(f, g_1,\ldots,g_n\right)(\xi).
\end{align}
This  is an immediate consequence of the fact that the operators $A_t$ and  $P_r$  commute, since both are  zonal Fourier  multipliers, and  of the Minkowski's inequality.  
\end{itemize}

\emph{Proof of Theorem \ref{pepe}}. The  necessary condition was already proved in Proposition \ref{teorema1}, taking $g_k=T_k((-\Delta_S)^kf)$, $k=1,2\ldots,n$. 

Now we prove the sufficient condition. We start with the case $2n < \alpha < 2(n+1)$. Suppose that there exists  $g_1,g_2, \ldots, g_n \in L^2 \left( \mathbb{S}^{d-1}  \right) $ such that $f \; \textrm{and } S_\alpha (f,g_1, g_2, \cdot \cdot \cdot, g_n ) \in L^2 \left( \mathbb{S}^{d-1}  \right) $.

For $r \in (0,1)$ and $\xi \in \mathbb{S}^{d-1}$,  we claim  that
	\begin{equation}\label{clain_1}
	P_rg_k(\xi)=T_k((-\Delta_S)^kP_rf)(\xi) ,\quad k=1, 2, \cdots, n,
	\end{equation}
	 which will be proved  later.

For $k=1,2, \cdot \cdot \cdot, n$ and $r\rightarrow 1^{-}$,  from ii) we have that $P_rg_k $ converges   to $g_k$ in $L^2 \left(  \mathbb{S}^{d-1}\right)$     and it is not difficult to check that $T_k((-\Delta_S)^kP_rf)$  converges   to $T_k((-\Delta_S)^kf)$ in  $\mathcal{D}'(  \mathbb{S}^{d-1} )$.  Then we conclude that $$ T_k((-\Delta_S)^k f)=g_k \in L^2(\mathbb{S}^{d-1} ).$$
Since $T_k$ is invertible in $L^2(\mathbb{S}^{d-1})$, it follows that $f \in H^{2n}\left( \mathbb{S}^{d-1} \right)$ and Theorem \ref{pepe} follows now from Proposition \ref{teorema1}.

 To prove (\ref{clain_1}), we introduce the auxiliary function
	\begin{align}\label{op:D} 
	D_\alpha (\eta)&:=\bigg(\int_0^{\pi} \bigg|
	\sum_{k=1}^n\,{T_k}((-\Delta_S)^k P_rf)(  \eta  ) A_t \big( \left|  \eta - \cdot \right|^{2k}    \big)( \eta )\nonumber\\
	&-\sum_{k=1}^nP_rg_k(  \eta)\, A_t \big( \left| \eta - \cdot \right|^{2k}    \big)( \eta )\bigg|^2\frac{dt}{t^{2\alpha+1}}\bigg)^{\frac{1}{2}}, \hspace{0.3cm} \eta \in \mathbb{S}^{d-1}.
	\end{align}
We have that $D_\alpha$ in $L^2(\mathbb{S}^{d-1})$, since from	 H\"older's inequaility and \eqref{propclave} we have
	\begin{align*}
	D_\alpha( \eta  )&\leq S_\alpha (P_rf,  T_1  ((-\Delta_S)P_rf),\ldots ,  T_n   ((-\Delta_S)^nP_rf))( \eta  ) \\ &+  S_\alpha (P_rf,P_rg_1, \ldots ,P_rg_n)(  \eta  )\\
	&\leq S_\alpha (P_rf,T_1((-\Delta_S)P_rf), \ldots ,T_n\left((-\Delta_S)^nP_rf\right)( \eta  ) \\ &+ P_r S_\alpha (f,g_1,  \ldots ,g_n)( \eta  ),
	\end{align*} 
 and $S_\alpha (P_rf,T_1((-\Delta_S))P_rf, \ldots, T_n ((-\Delta_S)^nP_rf)\in L^2(\mathbb{S}^{d-1})$ (is a consequence of  Proposition \ref{teorema1} since $P_rf \in H^\alpha \left(   \mathbb{S}^{d-1}  \right)$) and $S_\alpha (f,g_1,  \ldots ,g_n)$ is also in $L^{2} \left( \mathbb{S}^{d-1}  \right)$ (by hypothesis).

Without loss of generality we may assume that $D_\alpha (\xi)< \infty$. We are going to prove that if \eqref{clain_1} is false, then $D_\alpha (\xi)$ is not finite and we have a contradiction. To clarify the proof, we start by proving the case $k=1$ in (\ref{clain_1}) and the general case $1<k \leq  n$ is deduced by induction.
	
	 By (\ref{Atfunct}),
	\begin{equation}\label{zapato}
	A_t\left( \left| \xi - \cdot  \right|^{2k}   \right)(\xi) \backsim_{d,n} (1-\cos t)^{k}, \;\;\; k=1,2, \ldots , n,   
	\end{equation}
	then we have that
	$$D^2_\alpha (\xi)\sim_{d,n}\int_0^\pi \left| P_rg_1(\xi) -T_1  ((-\Delta_S))P_rf)(\xi)+G_1(\xi , t)     \right|^2 \frac{(1-\cos t)^2}{t^{2 \alpha +1 }}dt,$$
where
	 $$ \vert G_1(\xi,t)\vert\lesssim_{d,n}\sum_{j=2}^n
		\left| P_rg_j (\xi)-   T_j  ((-\Delta_S)^j P_rf)(\xi)  \right| (1-\cos t)^{j-1}.$$
	Hence, $\lim_{t\rightarrow 0+}G_1(\xi,t)=0$ so that  $P_rg_1(\xi) -  T_1  ((-\Delta_S)P_rf)(\xi)+G_1(\xi , t) $  can be extended in the variable $t$ to a continuous function on   $\left[0,\pi \right]$  with 
	$$\lim_{t\longrightarrow 0^+}\left( P_rg_1(\xi) - T_1  ((-\Delta_S)P_rf)(\xi)+G_1(\xi , t) \right)= P_rg_1(\xi) -   T_1  ((-\Delta_S)P_rf)(\xi).$$
	If $ \left| P_rg_1(\xi) -  T_1  ((-\Delta_S)P_rf)(\xi)  \right| \neq 0$, there exists $\kappa\in \left( 0, \pi  \right)$ such that if $t \in (0, \kappa)$ we have that 
	$$\left| P_rg_1(\xi) -   T_1  ((-\Delta_S)P_rf)(\xi)+G_1(\xi , t)   \right|^2 \geq \frac{1}{4} \left| P_rg_1(\xi) -   {T_1}  ((-\Delta_S)P_rf)(\xi)\right|^2,$$ and 
	\begin{equation}\label{contradiccion}
	\infty > D_\alpha^2(\xi)  \gtrsim_{d,n} \frac{1}{4} \left| P_rg_1(\xi) -  T_1  ((-\Delta_S)P_rf)(\xi)\right|^2 \int_0^{\kappa} \frac{(1-\cos t)^2}{t^{2 \alpha +1 }}dt.
	\end{equation}
	However, the last integral term in \eqref{contradiccion} is finite if and only if $\alpha <2$, which gives a contradiction. Therefore,
	$$P_rg_1(\xi) = T_1  ((-\Delta_S)P_rf)(\xi)  .$$
	
	Consider now the general case. Suppose that we have shown that 
	$$P_rg_j (\xi)= T_j  ((-\Delta_S)^jP_rf)(\xi)  , \;\; j=1,2, \ldots, k<n.$$
	To prove that $P_rg_{k+1}(\xi) =  {T_{k+1}}  ((-\Delta_S)^{k+1}P_rf)(\xi)$, we write 
	$$D_\alpha^2(\xi)\sim_{d,n}\int_0^\pi  \left| P_rg_{k+1}(\xi) - T_{k+1}  ((-\Delta_S)^{k+1}P_rf)(\xi)+G_k(\xi , t)     \right|^2 \frac{(1-\cos t)^{2k}}{t^{2 \alpha +1 }}dt,$$
	where
	$$ \vert G_k(\xi,t)\vert\lesssim_{d,n}\sum_{j=k+2}^n
	\left|  P_rg_j (\xi)-{T_j}((-\Delta_S)^j P_rf)(\xi)   \right|(1-\cos t)^{j-1}.$$
	
	If we repeat the process carried out for $k=1$ in \eqref{contradiccion}, we obtain that if $D_\alpha (\xi)< \infty$, then  $\alpha$  has to be less that $2k$, but this is not possible since $\alpha \geq 2n > 2k$. Therefore, we have proved (\ref{clain_1}).

	%Now, for $k=1,2, \cdot \cdot \cdot, n$ and $r\rightarrow 0+$  we have that $P_rg_k $ converges   to $g_k$ in $L^2 \left(  \mathbb{S}^{d-1}\right)$     and $c_k((-\Delta_S)^kP_rf)$  converges   to $c_k((-\Delta_S)^kf)$ in  $\mathcal{D}'\left(  \mathbb{S}^{d-1} \right)$, then we conclude that $$ c_k((-\Delta_S)^k_S f)=g_k \in L^2\left(\mathbb{S}^{d-1} \right)$$ for $k=1,2, \cdots n,$ and  since $c_k$ is invertible in $L^2(\mathbb{S}^{d-1})$ it follows that $f \in H^{2n}\left( \mathbb{S}^{d-1} \right)$. 
	
	%Theorem \ref{pepe} follows now from Proposition \ref{teorema1}.

	The remaining case $\alpha = 2n$ follows in a similar way using  the  auxiliary function in $\mathbb{S}^{d-1}$
	\begin{align}
	D_{2n}^2 (\eta)&=\int_0^{\pi} \bigg|\sum_{k=1}^{n-1}\, T_k   ((-\Delta_S)^k P_rf)( \eta   ) A_t \big( \left|  \eta - \cdot \right|^{2k} \big)(  {\eta} ) \nonumber \\ 
	&+  {T_n}    \left(A_t\left( (-\Delta_S)^n P_rf\right)\right)(  {\eta}  ) A_t \big( \left|  {\eta} - \cdot \right|^{2n}    \big)(  {\eta}  ) \nonumber  \\
	&-\sum_{k=1}^{n-1}P_rg_k(  {\eta}  )\, A_t \big( \left|  {\eta} - \cdot \right|^{2k} \big)(  {\eta} )-A_t \left(P_rg_n\right)(  {\eta}  )  A_t \big( \left|  {\eta} - \cdot \right|^{2n}   \big)(  {\eta}  )\bigg|^2\frac{dt}{t^{4n+1}}  \nonumber \\
	&=\int_0^{\pi} \bigg|\sum_{k=1}^{n-1}  \left(  {T_k}   ((-\Delta_S)^k P_rf)(  {\eta}  ) -P_rg_k(  {\eta}  )\right)   A_t \big( \left|  {\eta} - \cdot \right|^{2k} \big)(  {\eta} ) \nonumber   \\ \label{alegre}
	&+ \left( {T_n}   \left(A_t \left((-\Delta_S)^nP_rf\right)\right)( {\eta}  )
	- A_t \left(P_rg_n\right)(  {\eta}  ) \right)     A_t \big( \left|  {\eta} - \cdot \right|^{2n} \big)(   {\eta} ) \bigg|^2 \frac{dt}{t^{4n+1}}, %\label{alegre}
	\end{align} 
	that is a function in $L^2(\mathbb{S}^{d-1})$. Since $P_rg_n$ and ${T_n}  ((-\Delta_S)^nP_rf) $ are continuous, $T_n$ is a zonal Fourier multiplier, we have
	\begin{equation}\label{limite}
	\lim_{t\longrightarrow 0^+}\big( {T_n}    (A_t \left((-\Delta_S)^nP_rf)\right)(\xi)
	- A_t \left(P_rg_n\right)(\xi) \big)
	\end{equation}
	$$= {T_n}    ((-\Delta_S)^nP_rf)(\xi)  {-\left(P_rg_n\right)(\xi) },$$
	 for every $\xi \in \mathbb{S}^{d-1}$.
	
	We assume again that $D_{2n}(\xi)< \infty$. (\ref{limite})  together with the same argument that we applied to prove (\ref{clain_1}) would imply that 
	\begin{equation}\label{alegre_1}
	  {T_k}  ((-\Delta_S)^kP_rf)(\xi) =P_rg_k(\xi), \hspace{0.3cm} k=1,2, \ldots, n-1.
	  \end{equation}
	
	From (\ref{alegre}), (\ref{alegre_1}) and (\ref{Atfunct}) we have	
	$$\int_0^{\pi}\left|  {T_n}   (A_t \left((-\Delta_S)^nP_rf\right)(\xi)
	- A_t \left(P_rg_n\right)(\xi)      \big)(\xi) \right|^2 \frac{(1-\cos t)^{{2n}}}{t^{4n+1}}dt < \infty, $$ but for (\ref{limite}) this is possible only if 
$$  {T_n}   ((-\Delta_S)^nP_rf)(\xi) =P_rg_n(\xi),$$ since the integral $\int_0^{{\pi}}\frac{(1-\cos t)^{{2n}}}{t^{4n+1}}dt$  is not convergent.

Now proceeding as in the case $2n< \alpha < 2(n+1)$, we conclude that $f \in H^{2n} \left( \mathbb{S}^{d-1}  \right)$.
\hfill$\Box$

\section{Proof of auxiliary results}\label{auxiliar1}

\subsection{Proof of Lemma \ref{lemma:At} and Corollary \ref{corolarioAt}.} Since the proof of Lemma \ref{lemma:At} requires the use of the spherical harmonics, $Y_{\ell}\in \mathbb{H}_\ell^d$, we refer the interested reader to \cite{AH,EW,SW} to recall its main properties . In particular, we need two elements: the following representation
\begin{equation}\label{representaY}
Y_{\ell}(\xi):=\sum_{k=1}^{\nu(\ell)}a_k\,P_{\ell,d}(\eta_{k}\cdot\xi),
\end{equation}
where $a_k\in \mathbb{C}$ and $\eta_k\in\mathbb{S}^{d-1}$, and the next  lemma.
\begin{lemma}\label{unicidadHl}  Let $\xi\in\mathbb{S}^{d-1}$ and $L_\ell \in\mathbb{H}_\ell^d$ such that
 $L_\ell(R\eta)=L_\ell(\eta)$ for all rotations $R$ in $\mathbb{R}^d$ such that $R(\xi)=\xi$.
Then $L(\eta)=L(\xi)P_{\ell,d}(\eta\cdot\xi).$
\end{lemma}
The proof of this lemma follows the same lines as \cite[Theorem 4.10 ]{EW}.

\emph{Proof of Lemma \ref{lemma:At}.}
 To prove that $A_t$ is a zonal Fourier multiplier it is enough to see that for any real $\ell \geq 0$ and $Y_\ell\in \mathbb{H}_\ell^d$ we have
\begin{equation}\label{At:inter}
A_t Y_\ell=m_{\ell,t}Y_\ell,
\end{equation}
for $m_{\ell,t}$ defined in \eqref{mutliplierAt}. Fix $\xi\in\mathbb{S}^{d-1}$. Then, from (\ref{def:At}) and (\ref{representaY})
\begin{align}\label{representacion_1}
A_tY_\ell(\xi)&=\frac{1}{|C(\xi,t)|}\int_{C(\xi,t)}Y_\ell (\tau) d \sigma (\tau)\nonumber\\
&=\frac{1}{|C(\xi,t)|}\sum_{k=1}^{\nu (\ell)}a_k \int_{C(\xi,t)} P_{\ell , d}(\eta_k \cdot \tau) d \sigma (\tau) =\frac{1}{|C(\xi,t)|}\sum_{k=1}^{\nu (\ell)}a_k L(\eta_k),
\end{align}

where 
\begin{equation*}
L(\eta):=\int_{C(\xi,t)}P_{\ell,d}(\eta\cdot\tau)\,d\sigma(\tau).
\end{equation*}
We are going to check that $L$ satisfies the hypotheses of the Lemma \ref{unicidadHl}. From the addition theorem (see \cite[Theorem 2.9]{AH})
\[\sum_{j=0}^{\nu(\ell)}\overline{Y_{\ell}^j(\eta)}Y_\ell^j(\tau)=\frac{\nu(\ell)}{\vert\mathbb{S}^{d-1}\vert}P_{\ell,d}(\eta\cdot\tau)\]
it is clear that $L\in\mathbb{H}^\ell_d$. Moreover, if $R$ is a rotation  such that $R(\xi)=\xi$, we have that in particular $R$ leaves $C(\xi,t)$  invariant and hence 
\begin{align*}
L(R\eta)=&\int_{C(\xi,t)}P_{\ell,d}(R\eta\cdot\tau)\,d\sigma(\tau)=\int_{C(\xi, t)}P_{\ell, d}( \eta \cdot R^{-1} \tau) d \sigma (\tau)\\
&=\int_{C(\xi, t)}P_{\ell, d}( \eta \cdot  \tau) d \sigma (\tau)=L(\eta).
\end{align*}
If we apply Lemma \ref{unicidadHl}, (\ref{representacion_1}) and (\ref{representaY}), we get
$$A_tY_\ell (\xi)=\frac{L(\xi)}{|C(\xi,t)|}\sum_{k=1}^{\nu (\ell)}a_k P_{\ell , d}(\eta_k \cdot \xi)=\frac{L(\xi)}{|C(\xi,t)|}Y_\ell (\xi),$$
and thus to finish the proof it is enough to prove that $\frac{L(\xi)}{|C(\xi,t)|}=m_{\ell, t}$. To calculate $L(\xi)$ we need  to evaluate an  integral over $C(\xi,t)$. To do this, we define $s\in (\cos t, 1)$ such that $\xi\cdot\tau=s$ and we first integrate over the parallel $L_s=\{\tau\in\mathbb{S}^{d-1}:\tau\cdot\xi=s\}$, orthogonal to $\xi$. Then we obtain a function of  $s$ and we integrate over the interval $(\cos t, 1)$; precisely, we have that
\begin{align}\label{change_variables}
L(\xi)&=\int_{\cos t}^1\int_{\sqrt{1-s^2}\,\mathbb{S}^{d-2}}P_{\ell,d}(s)\,d\sigma(\theta)\frac{ds}{\sqrt{1-s^2}}\\
&=|\mathbb{S}^{d-2}|\int_{\cos t}^1P_{\ell,d}(s)\,(1-s^2)^{\frac{d-3}{2}}\,ds=|C(\xi,t)|m_{\ell, t}.\nonumber
\end{align}
The details of this method of integration are described for example in \cite[Appendix D2]{G}.
\hfill$\Box$

\emph{Proof of Corollary \ref{corolarioAt}.}  To prove  \eqref{medida_casquete} and \eqref{Cte}, we use the following estimate
\begin{align}\label{suela}
& \int_{\cos t}^1 (1-s)^k(1-s^2)^{\frac{d-3}{2}}ds \backsim \int_{\cos t}^1 (1-s)^{k +\frac{d-3}{2}}ds \nonumber \\ & =\frac{2(1-\cos t)^{k+\frac{d-1}{2}}}{2k+d-1}, \;\;t \in [0,\pi], \;\; k=0,1, \ldots ,
\end{align}
that it is clear for $t\in \left[0, \frac{\pi}{2}   \right]$ and for $t\in \left[ \frac{\pi}{2} ,  \pi  \right]$ it is a consequence of the fact that the function 
$$f(t)=\left\{  \begin{array}{ll}
\frac{ \int_{\cos t}^0 (1-s)^{k +\frac{d-3}{2}}ds  }{ \int_{\cos t}^0 (1-s)^k(1-s^2)^{\frac{d-3}{2}}ds  }, \hspace{0.2cm} t \in \left(  \frac{\pi}{2},\pi  \right] , \\ 1, \hspace{4.2cm} t=\frac{\pi}{2},
\end{array}     \right. $$ 
is positive and continuous in $\left[\frac{\pi}{2} ,\pi   \right]$. Then, using the definition given in \eqref{def:casquete}, (\ref{suela}) and evaluating $|C(\xi,t)|$ as we did in \eqref{change_variables}, we trivially obtain
\begin{equation*}
|C(\xi,t)|=|\mathbb{S}^{d-2}|\,\int_{\cos t}^1(1-s^2)^{\frac{d-3}{2}}\,ds\sim_d\int_{\cos t}^1(1-s)^{\frac{d-3}{2}}\sim_d t^{d-1},
\end{equation*}
which also proves \eqref{Cte}.

Equality \eqref{Atlap} is an immediate consequence of \eqref{At:inter} together with \eqref{laplaciano}.

Finally, to prove \eqref{Atfunct}, the rotational invariance allows to consider that $\xi$ has as coordinates $(0,\ldots,0,1)$  and then
\begin{align*}
A_t(|\xi-\cdot|^{2k})(\xi)&=\frac{2^k}{|C(\xi,t)|}\int_{C(\xi,t)}(1-\tau_d)^k\,d\sigma(\tau),
\end{align*}
where $\tau_d$ is the last coordinate of $\tau.$	Evaluating this integral as we explained in \eqref{change_variables} and using (\ref{suela}), we obtain the desired equality and estimate.\hfill$\Box$

\subsection{  Proof of Lemma \ref{prop:estima} and Lemma \ref{prop:estimaalpha2n}} \hfill

\emph{Proof of Lemma \ref{prop:estima}}. Let $\ell\geq 1$. We start by proving that $ I_{\alpha, n}(\ell) \lesssim_{d,n} \ell^{2 \alpha}$. We split the integral term in two pieces
\begin{equation}\label{estimacion_0}
I_{\alpha, n}(\ell) =\int_0^{a/\ell}|M_{\ell,t}|^2\,\frac{dt}{t^{2\alpha+1}}+\int_{a/\ell}^\pi|M_{\ell,t}|^2\,\frac{dt}{t^{2\alpha+1}}=I_{\alpha, n}^1(\ell)+I_{\alpha, n}^2(\ell),
\end{equation}
where $a>0$ is an absolute constant. In order to control the second term, we need to consider $a>1$, for example, $a=2$. 

To estimate $I_{\alpha, n}^1(\ell)$, using the Taylor expansion of the Legendre polynomials we have
\begin{equation}\label{taylorexpan}
P_{\ell,d}(s)=\sum_{k=0}^{n}\frac{(-1)^{k}P_{\ell,d}^{(k)}(1)}{k!}(1-s)^{k}+\frac{P_{\ell,d}^{(n+1)}(\tau(s))}{(n+1)!}(1-s)^{n+1},
\end{equation}
for some $\tau(s)\in(s,1)$  and  
\begin{equation}\label{cota_deriva}
\left| P_{\ell,d}^{(k)} (t) \right| \leq  P_{\ell,d}^{(k)} (1)    \backsim_{d,k}   \ell^{2k}, \hspace{0.3cm} k=0,1,\ldots,
\end{equation}
(see \cite[pp.58]{AH} and (\ref{cotaderaiva_0})). Then having into account definition \eqref{multiplicador1}, (\ref{suela}) and \eqref{cota_deriva}, we deduce that
\begin{align}
I_{\alpha, n}^1(\ell)& =\int_0^{\frac{ 2 }{\ell}}\bigg|C_{t,d}\int_{\cos t}^1\frac{P_{\ell,d}^{(n+1)}(\tau(s))}{(n+1)!}(1-s)^{n+1}(1-s^2)^{\frac{d-3}{2}}\,ds\bigg|^2\frac{dt}{t^{2\alpha+1}}\nonumber\\
&\lesssim_{d,n}\ell^{4n+4} \int_0^{\frac{  2 }{\ell}}\bigg|C_{t,d}\int_{\cos t}^1(1-s)^{n+1+\frac{d-3}{2}}\,ds\bigg|^2\frac{dt}{t^{2\alpha+1}}\nonumber\\
&\lesssim_{d,n}\ell^{4n+4}\int_0^{\frac{2} {\ell}}(1-\cos t)^{2n+2}\frac{dt}{t^{2\alpha+1}}.\label{inter:I1}
\end{align}
Taking the bound
\begin{align}\label{cota:cos}
 1- \cos t \leq  \frac{t^2}{2},\quad t\in( 0,\pi),
\end{align}
in \eqref{inter:I1} we conclude that
\begin{align}\label{cotaI1}
I_{\alpha, n}^1(\ell)\lesssim_{d,n}\ell^{4n+4}\int_0^{\frac{2}{\ell}}\frac{dt}{t^{2\alpha-4n-3}}\lesssim_{d,n} \ell^{2\alpha}, \;\; \; \textrm{whenever $\alpha<2(n+1)$}.
\end{align}
%whenever $\alpha<2(n+1)$.  

Now we study $I_{\alpha, n}^2(\ell)$. By the  mean value theorem, \eqref{cota_deriva},  \eqref{suela}, \eqref{cota:cos} and the fact that $t\ell>2$,  we have that
\begin{align*}
|M_{\ell,t}|&\lesssim_dC_{t,d}\int_{\cos t}^1\bigg[\ell^2(1-s)+\sum_{k=1}^{n}\frac{\ell^{2k}}{k!}(1-s)^{k}\bigg](1-s^2)^{\frac{d-3}{2}}\,ds\\
&\lesssim_d\sum_{k=1}^{n}\ell^{2k}t^{-d+1}\int_{\cos t}^1(1-s)^{k+\frac{d-3}{2}}\,ds\lesssim_d\sum_{k=1}^{n}\ell^{2k}t^{-d+1}(1-\cos t)^{k+\frac{d-1}{2}}\\
&\lesssim_d\sum_{k=1}^{n}\ell^{2k}t^{2k}\lesssim_d\ell^{2n}t^{2n}.
\end{align*}
Taking this expression in $I_{\alpha, n}^2(\ell)$, since we are assuming that $2n< \alpha$, we obtain
\begin{equation}\label{cotaI2}
I_{\alpha, n}^2(\ell) \lesssim_{d,n}\ell^{4n}\int_{\frac{2}{\ell}}^\pi\frac{t^{4n}}{t^{2\alpha+1}}\,dt \lesssim_{d,n} \ell^{2 \alpha}. 
\end{equation}
Putting together estimates \eqref{cotaI1} and \eqref{cotaI2} in \eqref{estimacion_0} we obtain the desired inequality. 

It remains to show that $I_{\alpha,n}(\ell)\gtrsim_{d,n} \ell^{2\alpha}$. 

From (\ref{multiplicador1}) and      \eqref{taylorexpan} we have that
\begin{equation}\label{teresa_1}
I_{\alpha, n}(\ell) \geq   \int_0^{a(\ell)} \bigg|C_{t,d} \int_{\cos t}^1 \frac{P_{\ell,d}^{(n+1)}(\tau (s))}{(n+1)!} (1-s)^{n+1}(1-s^2)^{\frac{d-3}{2}}ds \bigg|^2 \frac{dt}{t^{2 \alpha +1}},
\end{equation}  
for some $\tau (s) \in (s,1)$ and $a(\ell)<\pi/4$ a positive number to be chosen later.  The mean value theorem together with bound \eqref{cota_deriva} give
$$\big|   P_{\ell,d}^{(n+1)}(r)-P_{\ell,d}^{(n+1)}(1)  \Big| \leq \big|  P_{\ell,d}^{(n+2)}(\eta (r))  \big| (1-r) \leq P_{\ell,d}^{(n+2)}(1)(1-r),$$ for some $\eta (r) \in ( r,1)$. This implies that 
$$P_{\ell,d}^{(n+1)} (r) \geq   P_{\ell,d}^{(n+1)} (1) -P_{\ell,d}^{(n+2)} (1)(1-r).$$
If we take $r$ such that 
\[0 < 1-r \leq \frac{P_{\ell,d}^{(n+1)} (1) }{2 P_{\ell,d}^{(n+2)} (1) } \;  \Leftrightarrow \; r \in (1-k_{\ell,d,n},1)           ,\]
\begin{align}\label {constantek}
k_{\ell,d,n}=\frac{P_{\ell,d}^{(n+1)} (1) }{2 P_{\ell,d}^{(n+2)} (1) }=\frac{n+(d+1)/2}{(\ell +n+d+1)(\ell-n-1)}\sim_{d,n}\frac{1}{\ell^2}.
\end{align}
then  we have that 
\begin{align}\label{Juan1}
P_{\ell,d}^{(n+1)} (r) \geq P_{\ell,d}^{(n+1)} (1) /2.
\end{align}
%for $s\in\big(1-k_{\ell,d,n},1\big)$ with
%\begin{align}\label {constantek}
%k_{\ell,d,n}=\frac{P_{\ell,d}^{(n+1)} (1) }{2 P_{\ell,d}^{(n+2)} (1) }=\frac{n+d+1/2}{(\ell +n+d+1)(\ell-n-1)}\sim_{d,n}\frac{1}{\ell^2}.
%\end{align}
To estimate below $I_{\alpha , n}(\ell)$ we want to use (\ref{Juan1}) when $r=\tau (s)$, ($\tau (s)$ is the number in the integral  on (\ref{teresa_1})), so   we define  $a(\ell)$ as 
\begin{equation}\label{Juan_2}
\cos a (\ell)=1-k_{\ell,d,n}.
\end{equation}
With this choice we have
$$1-k_{\ell,d,n}= \cos a (\ell) < \cos t < s < \tau (s) <1 \; \Rightarrow \; \tau (s) \in (1-k_{\ell,d,n},1),$$
and $a(\ell)\sim_{d,n} 1/\ell$. 
Thus we have

%, with this choice we have

%Defining $a(\ell)$ as $\cos a (\ell)=1-k_{\ell,d,n}$ we can then ensure that the $\tau (s)$ in (\ref{teresa_1}) satisfies
%\begin{equation}\label{tau_s}
%P_\ell^{(n+1)}(\tau (s)) \geq \frac{P_\ell^{(n+1)} (1) }{2}.
%\end{equation}    
%Observe also that by definition of $a(\ell)$ and \eqref{constantek} we have that  $a(\ell)\sim_{d,n} 1/\ell$. Taking \eqref{tau_s}, \eqref{cota_deriva} and \eqref{cota:cos} in \eqref{teresa_1} we have
\begin{align*}
I_{\alpha, n}(\ell) &\gtrsim_{d,n} \ell^{4n+4}\int_0^{a(\ell)} \bigg|C_{t,d} \int_{\cos t}^1 (1-s)^{n+1}(1-s^2)^{\frac{d-3}{2}}ds \bigg|^2 \frac{dt}{t^{2 \alpha +1}}\\
&\gtrsim_{d,n} \ell^{4n+4}\int_0^{a(\ell)} (1-\cos t)^{2(n+1)}\frac{dt}{t^{2 \alpha +1}}\gtrsim_{d,n} \ell^{4n+4}\int_0^{a(\ell)}\frac{t^{4n+4}}{t^{2\alpha+1}}dt.
\end{align*}
Since $\alpha<2(n+1)$, we conclude that $I_{\alpha,n}(\ell)\gtrsim_{d,n}\ell^{2\alpha}$.\hfill$\Box$

\emph{Proof of Lemma \ref{prop:estimaalpha2n}.} Let $\ell\geq 1$. To show \eqref{clave3}, we prove first  $J_{n}(\ell) \lesssim_{d,n} \ell^{2 \alpha}$. We split into two parts as we did in Lemma \ref{prop:estima}.
\begin{align}\label{Jn_inter}
J_{n}(\ell) =\int_0^{1/\ell}|N_{\ell,t}|^2\,\frac{dt}{t^{4n+1}}+\int_{1/\ell}^{\pi}|N_{\ell,t}|^2\,\frac{dt}{t^{4n+1}}=J_{ n}^1(\ell)+J_{n}^2(\ell) .
\end{align}
For the term $J_{ n}^1(\ell)$, 
we use
\begin{align*}
&C_{t,d}\int_{\cos t}^1 P_{\ell,d }(s)(1-s^2)^{\frac{d-3}{2}}\,ds\\
&=C_{t,d}\int_{\cos t}^1 (P_{\ell,d} (s)-P_{\ell,d} (1))(1-s^2)^{\frac{d-3}{2}}\,ds +1
\end{align*}
(that it is a consequence of \eqref{agosto}),  together with \eqref{taylorexpan} and definitions \eqref{multiplicador2} and \eqref{multiplicador1}. We obtain,
\begin{align}\label{Jn1_inter}
J_{n}(\ell)^1\lesssim  J_{n}^{1,1} (\ell)+ I_{2n,n}^1(\ell),
\end{align}
where $ I_{2n,n}^1(\ell)$ was defined in \eqref{estimacion_0} (in this case we can take $a=1$) and 
\begin{align*}\label{J1}
J_{n}^{1,1} (\ell)&= \int_0^{1/\ell}\bigg|\frac{P_{\ell,d}^{(n)}(1)}{n!}\, C_{d,t}^2\int_{\cos t}^1 (P_{\ell,d}(s)-P_{\ell,d}(1)) (1-s^2)^{\frac{d-3}{2}}\,ds\nonumber\\
&\times \int_{\cos t}^1 (1-r)^n (1-r^2)^{\frac{d-3}{2}}\,dr  \bigg|^2 \frac{dt}{t^{4n+1}}.
\end{align*}
From \eqref{cotaI1},  we have that
\begin{equation}\label{cotaI2nn}
I_{2n,n}^1\lesssim_{d,n}\ell^{4n}.
\end{equation}
The term $J_{n}^{1,1}(\ell)$ follows applying the mean value theorem to $P_{\ell, d}$, namely, for some $\tau(s)\in(s,1)$, we have the estimate
\begin{align}\label{cotaJ2n1}
J_{n}^{1,1} (\ell)&\lesssim_{d,n}\ell^{4n} \int_0^{1/\ell}\bigg|C_{d,t}^2\int_{\cos t}^1 P_{\ell,d}^{(1)}(\tau(s))(s-1) (1-s^2)^{\frac{d-3}{2}}\,ds\nonumber\\
&\times \int_{\cos t}^1 (1-r)^n (1-r^2)^{\frac{d-3}{2}}\,dr  \bigg|^2 \frac{dt}{t^{4n+1}}\lesssim_{d,n}\ell^{4n+4}\int_0^{1/\ell} \frac{(1-\cos t)^{2+2n}}{t^{4n+1}}\,dt\nonumber\\
&\lesssim_{d,n}\ell^{4n+4}\int_0^{1/\ell} t^{3}\,dt\lesssim_{d,n}\ell^{4n}.
\end{align}
where we have also applied inequalities \eqref{cota_deriva}, \eqref{suela} and \eqref{cota:cos}. Taking estimates \eqref{cotaI2nn} and \eqref{cotaJ2n1} in \eqref{Jn1_inter} we get that
\begin{equation}\label{cotaJ2nn}
J_{n}^1(\ell)\lesssim_{d,n}\ell^{4n}.
\end{equation}
For the second term in \eqref{Jn_inter}, we can split as we did in \eqref{Jn1_inter}. Then
\begin{align}
J_{n}^2(\ell)&\lesssim_{d,n}\int_{\frac{1}{\ell}}^{\pi}\bigg|C_{t,d}\int_{\cos t}^1\bigg[P_{\ell,d}(s)-P_{\ell,d}(1)-\sum_{k=1}^{n-1}b_k(1-s)^{k}\bigg](1-s^2)^{\frac{d-3}{2}}\,ds\bigg|^2\frac{dt}{t^{4n+1}}\nonumber\\
&+\ell^{4n}\int_{\frac{1}{\ell}}^{\pi}\bigg|C_{t,d}^2\int_{\cos t}^1P_{\ell,d}(s)(1-s^2)^{\frac{d-3}{2}}\,ds \int_{\cos t}^1 (1-r)^n (1-r^2)^{\frac{d-3}{2}}\,dr  \bigg|^2 \frac{dt}{t^{4n+1}}\nonumber\\
&=I_{2n,n-1}(\ell)^2+J_{n}^{2,2}(\ell)
\label{Jn2_inter}
\end{align}
Since 	$2(n-1)<2n$, from \eqref{cotaI2} we obtain
\begin{equation}\label{raquel_5}
I^2_{2n,n-1}(\ell) \lesssim_{d,n}\ell^{4n}.
\end{equation} 
In order to estimate $J_{n}^{2,2}(\ell)$ we split the term two parts:
\begin{align}\label{inter:Jn22}
J_{n}^{2,2}(\ell)&=\ell^{4n}\int_{\frac{1}{\ell}}^{\pi/4}\bigg|C_{t,d}^2\int_{\cos t}^1P_{\ell,d}(s)(1-s^2)^{\frac{d-3}{2}}\,ds \int_{\cos t}^1 (1-r)^n (1-r^2)^{\frac{d-3}{2}}\,dr \bigg|^2 \frac{dt}{t^{4n+1}}\nonumber\\
&+\ell^{4n}\int_{\pi/4}^{\pi}\bigg|C_{t,d}^2\int_{\cos t}^1P_{\ell,d}(s)(1-s^2)^{\frac{d-3}{2}}\,ds \int_{\cos t}^1 (1-r)^n (1-r^2)^{\frac{d-3}{2}}\,dr  \bigg|^2 \frac{dt}{t^{4n+1}}\nonumber\\
&:=J_{n}^{2,2,1}(\ell)+J_{n}^{2,2,2}(\ell).
\end{align}
For the first term, we distinguish two cases: $d=2$ and $d\geq 3$. The case $d=2$ follows from the explicit form (see \cite[pp.38]{AH}) 
\begin{equation*}%\label{Chebychev}
P_{\ell,d}(t)=\cos(\ell\,\arccos t),\quad t\in[-1,1].
\end{equation*}
 Then taking this form in $J_{n}^{2,2,1}(\ell)$, and using estimates {(\ref{Cte})}, \eqref{suela} and \eqref{cota:cos} we obtain
\begin{align}\label{raquel_6d2}
J_{n}^{2,2,1}(\ell) &\lesssim \ell^{4n}\int_{\frac{1}{\ell}}^{\frac{\pi}{4}} \left|t^{-2}\int_{\cos t }^1\cos(\ell\arccos s)(1-s^2)^{-\frac{1}{2}}ds \int_{\cos t}^1 (1-r)^{n-\frac{1}{2}}\,dr   \right|^2 \frac{dt}{t^{4n+1}}\nonumber\\
&\lesssim_{d,n} \ell^{4n}\int_{\frac{1}{\ell}}^{\frac{\pi}{4}} \left|t^{-2}\,\frac{\sin \ell t}{\ell}\, t^{2n+1}  \right|^2 \frac{dt}{t^{4n+1}}\lesssim_{d,n}\ell^{4n-2}\int_{\frac{1}{\ell}}^{\frac{\pi}{4}}\frac{dt}{t^3}\lesssim_{d,n}\ell^{4n},
\end{align}
where we have used the change of variables $s=\cos(\theta/\ell)$ in the first integral term.
%{\color{red}{\bf Comentario:} de hecho, no se necesita  en este caso $d=2$ haberse quedado solo con el cachito entre $1/\ell$ y $\pi/4$ y funciona para estimar de manera completa $J_{n}(\ell)^{2,2}$ con $d=2$. No pensáis que quedar\'ia m\'as elegante que lo separemos desde el principio? Yo me inclino a pensar que s\'i, pero no opini\'on clara.}

For the case $d\geq 3$, we use the following assyptotic expansion for $P_{\ell,d}$, 
\begin{align}
P_{\ell,d}(\cos\theta)\sim_d\frac{1}{\pi^{\frac{1}{2}}\ell^{\frac{d-2}{2}}}\frac{1}{\sin^{\,\frac{d-2}{2}}\theta}\cos\Big((\ell+\frac{d-2}{2})\theta+(d-2)\frac{\pi}{4}\Big)+g(l),\label{asymp_legendre}
\end{align}
where $1/\ell  \lesssim  \theta \leq \pi/4$ and $|g(\ell)| \leq C \ell^{-\frac{d}{2}}$,
  with $C$ a absolute constant. This result is a generalization of the Laplace-Heine formula (see \cite[Theorem 8.21.8]{Sze}, since the polynomials  $P_{\ell,d}$ are proportional to the Jacobi polynomials).
To check that C is an absolute constant when $1/\ell  \lesssim  \theta \leq \pi/4$ see  \cite{FW}.

 From (\ref{asymp_legendre}) we have
$$\left| P_{\ell,d} (\cos \theta ) \right| (\sin\theta)^{d-2} \lesssim_d \frac{(\sin\theta)^{\frac{d-2}{2}}}{\ell^{\frac{d-2}{2}}}+\frac{(\sin\theta)^{d-2}}{\ell^{\frac{d}{2}}}\lesssim_d\frac{\theta^{\frac{d-2}{2}}}{\ell^{\frac{d-2}{2}}},\hspace{0.3cm} 1/\ell\lesssim \theta <\pi /4.$$
By using the change of variables $s=\cos\theta$ in $J_{n}^{2,2,1}(\ell)$ and applying the above inequality together with \eqref{suela}, \eqref{cota_deriva} and \eqref{cota:cos} we obtain
\begin{align}\label{raquel_6}
J_{n}^{2,2,1}(\ell) &\lesssim_{d,n} \ell^{4n}\int_{\frac{1}{\ell}}^{\frac{\pi}{4}} \left|  (1- \cos t)^{n-1-\frac{d-3}{2}}\int_{ 0}^t P_{\ell,d} ( \cos \theta) (\sin \theta)^{d-2} d \theta   \right|^2 \frac{dt}{t^{4n+1}}\nonumber\\
&\lesssim_{d,n}\ell^{4n}\int_{\frac{1}{\ell}}^{\frac{\pi}{4}}\left| t^{2n-d+1}\bigg(\int_{0}^{\frac{1}{\ell}}\theta^{d-2}\,d\theta+\int_{\frac{1}{\ell}}^t\ell^{\frac{2-d}{2}}\theta^{\frac{d-2}{2}}\,d\theta\bigg)\right|^2\frac{dt}{t^{4n+1}}\nonumber\\
&\lesssim_{d,n}\ell^{4n}\bigg(\frac{1}{\ell^{2d-2}}\int_{\frac{1}{\ell}}^{\frac{\pi}{4}}  \frac{dt}{t^{2d-1}} +\frac{1}{\ell^{d-2}}\int_{\frac{1}{\ell}}^{\frac{\pi}{4}}\frac{dt}{t^{d-1}} \bigg) 
\lesssim_{d,n}\ell^{4n}.
\end{align}

The second term in \eqref{inter:Jn22} is trivially bounded using \eqref{suela}, \eqref{cota_deriva} and \eqref{cota:cos} since 
\begin{align}\label{falta}
J_{n}^{2,2,2}\lesssim_{d,n}\ell^{4n}\int_{\frac{\pi}{4}}^\pi \left|C_{t,d}^2\,t^{d-1}\,t^{2n+d-1} \right|^2\frac{dt}{t^{4n+1}}\lesssim_{d,n}\ell^{4n}\int_{\frac{\pi}{4}}^\pi\frac{dt}{t}\lesssim_{d,n}\ell^{4n}.
\end{align}
Putting together \eqref{cotaJ2nn}, \eqref{Jn2_inter}, \eqref{raquel_5}, \eqref{inter:Jn22}, \eqref{raquel_6} (alternatively \eqref{raquel_6d2} in the case $d=2$) and \eqref{falta} in \eqref{Jn_inter}, we conclude that $J_{n}(\ell)     \lesssim_{d,n} \ell^{4n}$.

It remains to prove that $J_n(\ell) \gtrsim_{d,n} \ell^{4n}$.  Denote $h_n=A_t(|\xi-\cdot|^{2n})(\xi)$.  Using in \eqref{multiplicador2} the Taylor approximations of $P_{\ell , d}$ of order $n-1$ and $2$ at $s=1$, we obtain
\begin{align}
N_{\ell,t}&=C_{t,d}\int_{\cos t}^1\bigg(\frac{P_{\ell,d}^{(n)}(1)}{n!}(s-1)^n+\frac{P_{\ell,d}^{(n+1)}(\tau_1(s))}{(n+1)!}(s-1)^{n+1}\bigg)(1-s^2)^{\frac{d-3}{2}}ds\nonumber\\
&-\frac{(-1)^nP_{\ell,d}^{(n)}(1)}{n!}\,\frac{h_nC_{t,d}}{2^n}\int_{\cos t}^1P_{\ell,d}(1)(1-s^2)^{\frac{d-3}{2}}ds\nonumber\\
&-\frac{(-1)^nP_{\ell,d}^{(n)}(1)}{n!}\,\frac{h_n C_{t,d}}{2^n}\int_{\cos t}^1P_{\ell,d}^{(1)}(\tau_2(s))(s-1)(1-s^2)^{\frac{d-3}{2}}ds\nonumber\\
&=C_{t,d}\int_{\cos t}^1\frac{P_{\ell,d}^{(n+1)}(\tau_1(s))}{(n+1)!}(s-1)^{n+1}(1-s^2)^{\frac{d-3}{2}}ds\nonumber\\
&-\frac{(-1)^nP_{\ell,d}^{(n)}(1)}{n!}\,\frac{C_{t,d}h_n}{2^n}\int_{\cos t}^1P_{\ell,d}^{(1)}(\tau_2(s))(s-1)(1-s^2)^{\frac{d-3}{2}}ds  \nonumber\\
&:=N_{\ell,t}^1+N_{\ell,t}^2,\label{Ninterm}
\end{align}

where $\tau_1(s), \tau_2(s) \in (s,1)$ and we have used that $P_{\ell, d}(1)=1$ and 
$$\frac{P_{\ell,d}^{(n)}(1)}{n!}C_{t,d}\int_{\cot t}^1(s-1)^n(1-s^2)^{\frac{d-3}{2}}ds $$ $$=-\frac{(-1)^nP_{\ell,d}^{(n)}(1)}{n!}\,\frac{h_nC_{t,d}}{2^n}\int_{\cos t}^1(1-s^2)^{\frac{d-3}{2}}ds.$$

 Hence
\begin{align}\label{Jdosterminos}
J_n(\ell)^{\frac{1}{2}}\geq\bigg(\int_0^{c(\ell)}|N_{\ell,t}^2|^2\frac{dt}{t^{4n+1}}\bigg)^{\frac{1}{2}}-\bigg(\int_0^{c(\ell)}|N_{\ell,t}^1|^2\frac{dt}{t^{4n+1}}\bigg)^{\frac{1}{2}},
\end{align}
with $c(\ell)<\pi/2$ a positive number to be chosen later. The two terms of the right-hand side of the inequality above are  of the order  of $\ell^{4n}$ as $\ell  \rightarrow \infty$, but the first, as we will see,  will absorb the second one. To see this we will use a similar argument to that used in the proof of Lemma \ref{prop:estima}  to prove that $I_{\alpha ,n }(\ell)\gtrsim_{d,n} \ell^{2 \alpha}$. More precisely, using the mean value theorem together with \eqref{cota_deriva} we have that
\[P_{\ell,d}^{(1)}(r)\geq P_{\ell,d}^{(1)}(1)-(1-r)P_{\ell,d}^{(2)}(1).\]
Then if we take $r$ such that 
$$0 < 1-r \leq \frac{P_{\ell,d}^{(1)} (1)} {b P_{\ell,d}^{(2)} (1)},$$ 
where $b>1$   will be a fixed number that  will depend on $d$ and $n$ and will be chosen later, we can assure that the following condition is satisfied
\begin{equation}\label{ventana}
P_{\ell,d}^{(1)}(r)\geq \left(1-\frac{1}{b} \right)P_{\ell,d}^{(1)}(1)>0, \hspace{0.4cm} r \in \left[1-{P_{\ell,d}^{(1)} (1)} /({b P_{\ell,d}^{(2)} (1)}), 1 \right).
\end{equation}
To give a lower bound  for the first term in the right-hand side of \eqref{Jdosterminos}, we need condition (\ref{ventana}) to be satisfied when we take $r$  as the value  $\tau_2(s)$ that appears in $N_{\ell , t}^2$. To this aim we define  $c(\ell)$,  for $\ell \geq 3$, as
\begin{align}\label{escojoA}
&\cos c (\ell)=1-\frac{P_{\ell,d}^{(1)} (1)}{ b P_{\ell,d}^{(2)} (1) }=1- \frac{2 \left( 1+ \frac{d-1}{2}     \right)}{b(\ell -1)(\ell+d-1)} \nonumber   \\ &\geq \left( 1-\frac{d+1}{2bd} \right)  \geq \left( 1-\frac{1}{b}   \right),
\end{align}
then
$$1-\frac{P_{\ell,d}^{(1)} (1)}{ b P_{\ell,d}^{(2)} (1) }=\cos c(\ell) < \tau_2(s)<1 \; \Rightarrow \; \tau_2(s) \in  \left[1-{P_{\ell,d}^{(1)} (1)} /({b P_{\ell,d}^{(2)} (1)}), 1 \right),
$$ and from (\ref{escojoA})
\begin{equation}\label{c(ele)}
c(\ell)\backsim_{d, b} \frac{1}{\ell}.
\end{equation}
In order to obtain an  upper and a lower bound for $\left| N_{\ell , t}^1  \right|$ and $\left| N_{\ell , t}^2  \right|$ respectively, we use the following inequalities  valid whenever $s \in (\cos t, 1)$, $t \in (0, c (\ell))$ and $\ell \geq  2$. Since
\begin{equation}\label{uno} 
2^{\frac{d-3}{2}} \left(1-\frac{1}{b}\right)^{\frac{d-3}{2}} \leq (1+s)^{\frac{d-3}{2}} \leq 2^{\frac{d-3}{2}}, \hspace{0.3cm} d \geq 3,
\end{equation}
we have that
\begin{align}
&\frac{2^{\frac{d-1}{2}}}{2j+d-1}    \left(1-\frac{1}{b}\right)^{\frac{d-3}{2}}    (1- \cos t)^{j+\frac{d-1}{2}}    \leq   \int_{\cos t}^1 (1-s)^j(1-s^2)^{\frac{d-3}{2}}ds\nonumber\\
&\leq \frac{2^{\frac{d-1}{2}}}{2j+d-1}(1- \cos t)^{j+\frac{d-1}{2}}, \hspace{0.4cm} j=0,1, \ldots, \; \; d \geq 3.\label{dos}
\end{align}
Similarly, in the particular case $d=2$ we have 
\begin{align}\label{tres}
  \frac{(1-\cos t)^{j+\frac{1}{2}}}{\sqrt{2}\left(  j+\frac{1}{2}  \right)} &\leq   \int_{\cos t}^1 (1-s)^j(1-s^2)^{-\frac{1}{2}}ds\nonumber\\
  &\leq  \left( 1-\frac{1}{b}    \right)^{-\frac{1}{2}}  \frac{(1-\cos t)^{j+\frac{1}{2}}}{\sqrt{2}\left(  j+\frac{1}{2}  \right)} 
\end{align}
for $j=0,1, \ldots $.

We start estimating from below the first integral of (\ref{Jdosterminos}) for $d\geq 3$. From (\ref{medida_casquete}), the expression of $h_n$ given by (\ref{Atfunct}), (\ref{ventana}) for $r=\tau_2(s)$ and   (\ref{dos}), we obtain

\begin{align*}
\int_0^{c(\ell)} \left| N^2_{\ell , t} \right|^2 \frac{dt}{t^{4n+1}} &\geq \frac{\left(  P_{\ell , d}^{(n)}(1)   P_{\ell , d}^{(1)}(1)    \right)^2}{(n!)^2}\left(1-\frac{1}{b}\right)^2\nonumber\\
&\times \int_0^{c(\ell)}\left| \frac{ \int_{\cos t}^1 (1-v)^n(1-v^2)^{\frac{d-3}{2}}dv \int_{\cos t}^1 (1-u)(1-u^2)^{\frac{d-3}{2}}du      }{\left(\int_{\cos t}^1 (1-s^2)^{\frac{d-3}{2}}ds\right)^2}       \right|^2 \frac{dt}{t^{4n+1}}\nonumber\\
&\geq \frac{\left(  P_{\ell , d}^{(n)}(1)   P_{\ell , d}^{(1)}(1)    \right)^2(d-1)^4}{(n!)^2(2n+d-1)^2(d+1)^2}                           \left(1-\frac{1}{b}\right)^{2(d-2)} \int_0^{c(\ell)} \frac{(1- \cos t)^{2n+2}}{t^{4n+1}}dt.
\end{align*}
Since
$$1-\cos t\geq\frac{\cos c(\ell)}{2}t^2,$$
whenever $t\in(0,c(\ell))$, the preceding integral becomes
\begin{equation}\label{primera}
\int_0^{c(\ell)} \left| N^2_{\ell , t} \right|^2 \frac{dt}{t^{4n+1}}\geq    
\frac{\left(  P_{\ell , d}^{(n)}(1)   P_{\ell , d}^{(1)}(1)    \right)^2(d-1)^4 c(\ell)^4 \cos^{2n+2} c(\ell)}{(n!)^2(2n+d-1)^2(d+1)^22^{2n+2}}  
\left(1-\frac{1}{b}\right)^{2(d-2)}.
\end{equation}

Now we estimate from above the second integral of (\ref{Jdosterminos}). Using \eqref{cota_deriva}, (\ref{cota:cos}) and (\ref{dos}), we obtain
\begin{align}\label{segunda}
& \int_0^{c(\ell)} \left| N^1_{\ell , t} \right|^2 \frac{dt}{t^{4n+1}} \nonumber \\ & \leq \frac{\left(  P_{\ell , d}^{(n+1)}(1)  \right)^2}{((n+1)!)^2} \int_0^{\cos t} \left|C_{\ell , t} \int_{\cos t}^1 (1-s)^{n+1}(1-s^2)^{\frac{d-3}{2}}ds      \right|^2
\frac{dt}{t^{4n+1}}  \nonumber \\
&
 \leq \frac{\left(  P_{\ell , d}^{(n+1)}(1)  \right)^2(d-1)^2}{((n+1)!)^2(2n+d+1)^2} \left(1-\frac{1}{b}\right)^{-(d-3)}\int_0^{c(\ell)} \frac{(1- \cos t)^{2n+2}}{t^{{4}n+1}} dt 
\nonumber \\
 & \leq \frac{\left(  P_{\ell , d}^{(n+1)}(1)  \right)^2(d-1)^2c (\ell)^4}{((n+1)!)^2(2n+d+1)^22^{2n+4}} \left(1-\frac{1}{b}\right)^{-(d-3)}.
 \end{align}
 
 From (\ref{Jdosterminos}),  (\ref{primera}), (\ref{escojoA}), (\ref{segunda}) and (\ref{cotaderaiva_0})  can be seen that
\begin{align}\label{2_parte}
& J_n(\ell)^{\frac{1}{2}} \geq \frac{ \ell ! (\ell +n+d-3)! \ell (\ell +d-2)c(\ell)^2(d-1) \Gamma \left(  \frac{d-1}{2} \right)       }{(\ell -n)! (\ell +d-3)!2^{2n+3} n! (d+1)\Gamma \left( n+\frac{d+1}{2}  \right)     } \left(  1-\frac{1}{b}  \right)^{-\frac{d-3}{2}}  \nonumber \\
&
 \times \left(   \left(  1-\frac{1}{b}  \right)^{n+\frac{2d-5}{2}}  -\frac{ (\ell -n) ( \ell +n +d -2  )(d+1)   }{ \ell ( \ell +d-2   )(n+1)( 2n+d+1 )}      \right) .
 \end{align}
As 
$$\lim_{\ell \longrightarrow \infty} \frac{ (\ell -n) ( \ell +n +d -2  )(d+1)   }{ \ell ( \ell +d-2   )(n+1)( 2n+d+1 )} =\frac{d+1   }{ (n+1)( 2n+d+1 )} ,$$ and $ n \geq 1$, there exists a $\ell (d,n)$ such that  if $\ell \geq \ell (d,n) $ we have
\begin{equation}\label{eleccion_ele}
 \frac{ (\ell -n) ( \ell +n +d -2  )(d+1)   }{ \ell ( \ell +d-2   )(n+1)( 2n+d+1 )} \leq \frac{d+1   }{ (n+1)( 2n+d+1 )}+\frac{1}{4} \leq \frac{3}{4}.
\end{equation} 
Now we take $b$, depending only on $d$ and $n$, sufficiently close to 1 satisfying that 
\begin{equation}\label{eleccion_b}
 \left(  1-\frac{1}{b}  \right)^{n+\frac{2d-5}{2}}   \geq \frac{3}{4}+\frac{1}{8}.
\end{equation}
From (\ref{c(ele)})
$$\frac{ \ell ! (\ell +n+d-3)! \ell (\ell +d-2)c(\ell)^2(d-1) \Gamma \left(  \frac{d-1}{2} \right)       }{(\ell -n)! (\ell +d-3)!2^{2n+3} n! (d+1)\Gamma \left( n+\frac{d+1}{2}  \right)     } \left(  1-\frac{1}{b}  \right)^{-\frac{d-3}{2}} \backsim_{d,n} \ell^{2n},$$
then this estimate, (\ref{2_parte}), (\ref{eleccion_ele}) and (\ref{eleccion_b}) show that  $J_n(\ell) \gtrsim_{d,n} \ell^{4n}$ for $d \geq 3$.

It is not difficult to verify that if $d=2$ we have
\begin{align}\label{3_parte}
& J_n(\ell)^{\frac{1}{2}} \geq \frac{ \ell ! (\ell +n-1)! \ell^2c(\ell)^2\pi      }{(\ell -n)! (\ell -1)!2^{n+1} n! \Gamma \left( n+\frac{3}{2}  \right)     } \left(  1-\frac{1}{b}  \right)^{-1}  
\nonumber \\
& \times \left( \frac{2}{3}  \left(  1-\frac{1}{b}  \right)^{n+4} -\frac{ (\ell -n) ( \ell +n   )   }{ \ell^2 (n+1)\left(  2n+3 \right)}      \right) ,
 \end{align}
and if we continue the calculations carried out for $d\geq 3$ we get $J_n(\ell) \gtrsim_{d,n} \ell^{4n}$ for $d =2 $.
\hfill$\Box$

\paragraph{\bf Acknowledgements} The first author was supported by the Spanish Grant MTM2017-85934-C3-3-P.  The second by the Spanish Grant MTM2017-82160-C2-1-P. The third by the Mexican Grant DGAPA-UNAM PAPIIT IN106418.

	%-------------------------------------------------------------------------

\end{document}